\documentclass[12pt]{amsart}
\usepackage{amssymb,amsmath}
\usepackage{subfigure}
\usepackage{graphpap,latexsym,epsf}
\usepackage{color,graphicx,psfrag}

\newtheorem{rem}{Remark}[section]

\newtheorem{teo}{Theorem}[section]
\newtheorem{defin}{Definition}[section]

\newtheorem{coro}{Corollary}[section]
\newtheorem{lemma}{Lemma}[section]

\def\proof{\noindent{\bf Proof.}\ }
\def\endproof{{\mbox{}\nolinebreak\hfill\rule{2mm}{2mm}\par\medbreak} }

\def\neweq#1{\begin{equation}\label{#1}}
\def\endeq{\end{equation}}
\def\bbbr{\mathbb{R}}

\def\supp{\mathop{\rm supp}}

\def\l{\lambda}

\def\O{\Omega}

\def\dis{\displaystyle}

\def\eps{\varepsilon}

\def\R{\bbbr}

\def\dis{\displaystyle}
\def\deg{{\rm deg}}

\def\hatu{{\widehat u}}

\def\D{\Delta}

\def\UB{{\Omega^0}}
\def\varkappa{\kappa}
\def\vk{\kappa}

\newcommand{\Bi}{B_i}
\newcommand{\ov}{\overline}
\newcommand{\id}{\mathop{\rm Id}}

\title[Strongly Competing Species in Special Domains]{ Coexistence and
  Segregation for Strongly Competing Species in Special Domains.}

\author[M.\ Conti, V. \ Felli]
{Monica Conti, Veronica Felli}

\address{
Monica Conti
\newline\indent
Dipartimento di Matematica ``F.Brioschi''
\newline\indent
Politecnico di Milano
\newline\indent
Via Bonardi 9, I-20133 Milano, Italy
\newline\indent
\emph{monica.conti@mate.polimi.it},}
\address{
Veronica Felli
\newline\indent
Dipartimento di Matematica e Applicazioni
\newline\indent
Universit\`a degli Studi di Milano-Bicocca
\newline\indent
Via Cozzi 53, I-20125 Milano, Italy
\newline\indent
\emph{veronica.felli@unimib.it}.}

\thanks{Research partially supported
the Italian MIUR Research Projects
{\it Metodi Variazionali ed Equazioni Differenziali Nonlineari},
}
\date{Revised version of November 16, 2007}

\begin{document}
\begin{abstract}
  We deal with strongly competing multispecies systems of
  Lotka-Vol--terra type with homogeneous Dirichlet boundary
  conditions. For a class of nonconvex domains composed by balls connected with
  thin corridors, we show the occurrence of pattern
  formation (coexistence and spatial segregation of all the species),
  as the competition grows indefinitely. As a result we prove the
  existence and uniqueness of solutions for a remarkable system of
  differential inequalities involved in segregation phenomena and
  optimal partition problems.

\end{abstract}

\maketitle
\section{Introduction}
In this paper we consider  the system of
$k\geq 2$ elliptic equations
\begin{equation}\label{modelLK}
%\left\{
%\begin{array}{l}
-\Delta u_i= \displaystyle f_i(x,u_i)-\vk u_i\sum_{j\neq i}
u_j\qquad\text{ in } \Omega,
%\end{array}\right.
\end{equation}
for $i=1,\dots,k$, where $\Omega\subset\R^N$ is a smooth, connected,
bounded domain.  Systems of this form model the steady states of $k$
organisms which coexist in the area $\Omega$.  The function $u_i$
represents the population density of the $i$-th species (hence only
$u_i\geq 0$ are considered) and $f_i$ describes the internal dynamic
of $u_i$.  The coupling between different equations is the classical
Lotka-Volterra interaction term: the positive constant $\varkappa$
prescribes the competitive character of the relationship between $u_i$
and $u_j$ and its largeness measure the strength of the competition.

Systems of this form have attracted considerable attention both in
ecology and social science since they furnish a relatively simple
model to study phenomena of extinction, coexistence and segregation
states of populations.  Several theoretical studies have been carried
out in this direction,
mainly in the case of two competing species and for the logistic
nonlinearities $f_i(u)=u(a_i-u)$. We quote for instance \cite{eflg,gl,kl,lm,mhmm,skt},
where it is shown that both coexistence and exclusion may occur, depending on the relations
between the diffusion rates, the coefficients of intra--specific and
 of inter--specific competitions.

In this paper we face the multispecies Lotka-Volterra system
(\ref{modelLK}) in the different perspective investigated by
\cite{ctv-asymp,dd3,dg1,dhmp} and we study the possibility
of coexistence governed by very strong competition.  As we shall
discuss in detail in Section~\ref{asintotic}, the presence of large
interactions of competitive type produces the spatial segregation of
the densities in the limit configuration as $\kappa\to \infty$, namely
if $(u_i^{\kappa})_{i=1,\dots,k}$ solves~(\ref{modelLK}), then, for all
$i=1,\dots,k$, $u_i^{\kappa}$ converges to some $u_i$ in $H^1(\Omega)$
which satisfies
\begin{equation}\label{segrego}
  u_i(x)\cdot u_j(x)=0\text{ a.e. in }\Omega,\qquad\text{for all } i\neq j,
\end{equation}
so that $\{u_i>0\}\cap\{u_j>0\}=\emptyset$.
Furthermore, in the limit, the densities
satisfy a system of differential inequalities of the form
\begin{equation}\label{2kvar}\left\{
    \begin{array}{l}
      -\Delta u_i\leq f_i(x,u_i),\qquad\text{ in } \Omega,\\
      -\Delta \widehat u_i\geq
      \widehat f_i(x,\widehat u_i),\qquad\text{ in } \Omega,\\
    \end{array}\right.
\end{equation}
where $\widehat{u}_i:=u_i-\sum_{j\neq i}u_j$ and $\widehat
f_i(x,\widehat u_i):=f_i(x,u_i)-\sum_{j\neq i} f_j(x,u_j)$, in the
sense of definition \ref{def:dise}.
The link between the differential inequalities \eqref{2kvar} and
population dynamics is reinforced by considering another class of
segregation states between species, governed by a minimization
principle rather than strong competition--diffusion.  In
\cite{ctv-var} (see also \cite{ctv2,ctv-fucik}), the following energy
functional
\begin{equation*}
  J(U)=\sum_{i=1,\dots, k}\left\{\int_\Omega\left(\frac12
|\nabla u_i(x)|^2- F_i(x,u_i(x))\right)dx\right\},
\end{equation*}
given by the sum of the internal energies of $k$ positive densities
$u_i$ having internal potentials $F_i(x,s)=\int_0^s f_i(x,u)du$, was
considered. The problem of finding the \emph{minimum} of $J(U)$ in the
class of $k$-tuples $U=(u_1,\dots,u_n)$ satisfying $u_j\cdot u_i=0$
a.e. on $\Omega$ for $i\neq j$ was investigated in \cite{ctv-var},
where it is proved that any non-trivial minimizer $U$ (if it exists)
satisfies the differential inequalities \eqref{2kvar}.

This further motivates the study of the solutions of
(\ref{segrego}--\ref{2kvar}) as a natural step in the understanding
of segregation phenomena occurring in population dynamics.  Remarkably
enough, \eqref{2kvar} coupled with \eqref{segrego} can be naturally
interpreted as a \emph{free boundary problem} with multiple phases: the unknown free
boundary set is given by
\[
{\mathcal F}=\bigcup_{i=1}^k\partial\{x\in\Omega:\,u_i(x)>0\},
\]
which represents the collection of the boundaries of the disjoint
supports of the densities. On its support each density $u_i$ solves
the elliptic equation $-\Delta u_i=f_i(x,u_i)$, while
the free boundary conditions are implicitly contained in
the global differential inequalities
\eqref{2kvar}.
The study of the properties of ${\mathcal F}$ is important from the
ecological point of view since it provides
information about how the segregation
occurs, in particular about the way the territory is partitioned by
the segregated populations.  In this direction, in
\cite{ctv,ctv2,ctv-var,ctv-asymp,ctv-uniq} a number of qualitative
properties both of $u_i$ and the free boundary set ${\mathcal F}$ is
exhibited.  We refer the interested reader to \cite{spt} for a brief
review of the regularity theory  so far developed,
and to the above quoted papers for proofs and details.

Another question of particular interest is the existence of a
\emph{strictly positive} solution to (\ref{segrego}--\ref{2kvar}),
that is a solution of the differential inequalities (\ref{2kvar})
with each
component $u_i\geq 0$ and $u_i$ positive on a set of positive measure.
As a matter of fact, since all the asymptotic states of the
Lotka-Volterra system have to satisfy (\ref{segrego}--\ref{2kvar}),
the existence of such a solution is necessary to ensure that
\emph{all} the species survive under strong competition.  It has to be
stressed that in \cite{ctv-var,ctv-asymp}, the strict positivity is
guaranteed by assuming positive boundary values for each component, in
the form $u_i=\phi_i\text{ on }\partial\Omega$ with $\phi_i>0$ on a set
of positive $(N-1)$--measure.

Hence a major problem consists in proving the existence of a positive
solution under natural boundary conditions, such as Dirichlet or
Neumann homogeneous boundary conditions.  This is precisely the
problem we face in this paper: we consider
(\ref{segrego}--\ref{2kvar}) with the Dirichlet condition
\begin{equation}\label{dirich}
u_i=0\qquad\text{ on }\partial\Omega,
\end{equation}
and we look for a strictly positive solution $U=(u_1,\dots,u_k)$.  The
interesting case of Neumann condition will be treated elsewhere (see
the concluding remarks).

This is an interesting and mathematically challenging problem: we
cannot expect in general to avoid extinction of one or more species.
For instance, if $\Omega$ is convex, it is shown
in \cite{kw} that two competing species can not coexist under strong
competition.  On the other end, the main variational procedure leading to
solutions of \eqref{2kvar}, that is the minimization of the internal
energy $J$, may fail under Dirichlet homogeneous conditions, since it
in general provides $k$-tuple of the form $(0,\dots, u_i,\dots 0)$
with all but one component identically zero, see \cite{sweers} for a
similar result.
Therefore, some mechanism of different nature must occur in order to
ensure coexistence of the species.

\begin{figure}[h]
  \centering \subfigure [the set $\O^0=B_1\cup B_2\cup B_3$ and
  segments $E$ joining the balls]{
\begin{psfrags}
\psfrag{E}{$E$}\psfrag{B1}{$B_1$}\psfrag{B2}{$B_2$}
\psfrag{B3}{$B_3$}
\includegraphics[height=4cm]{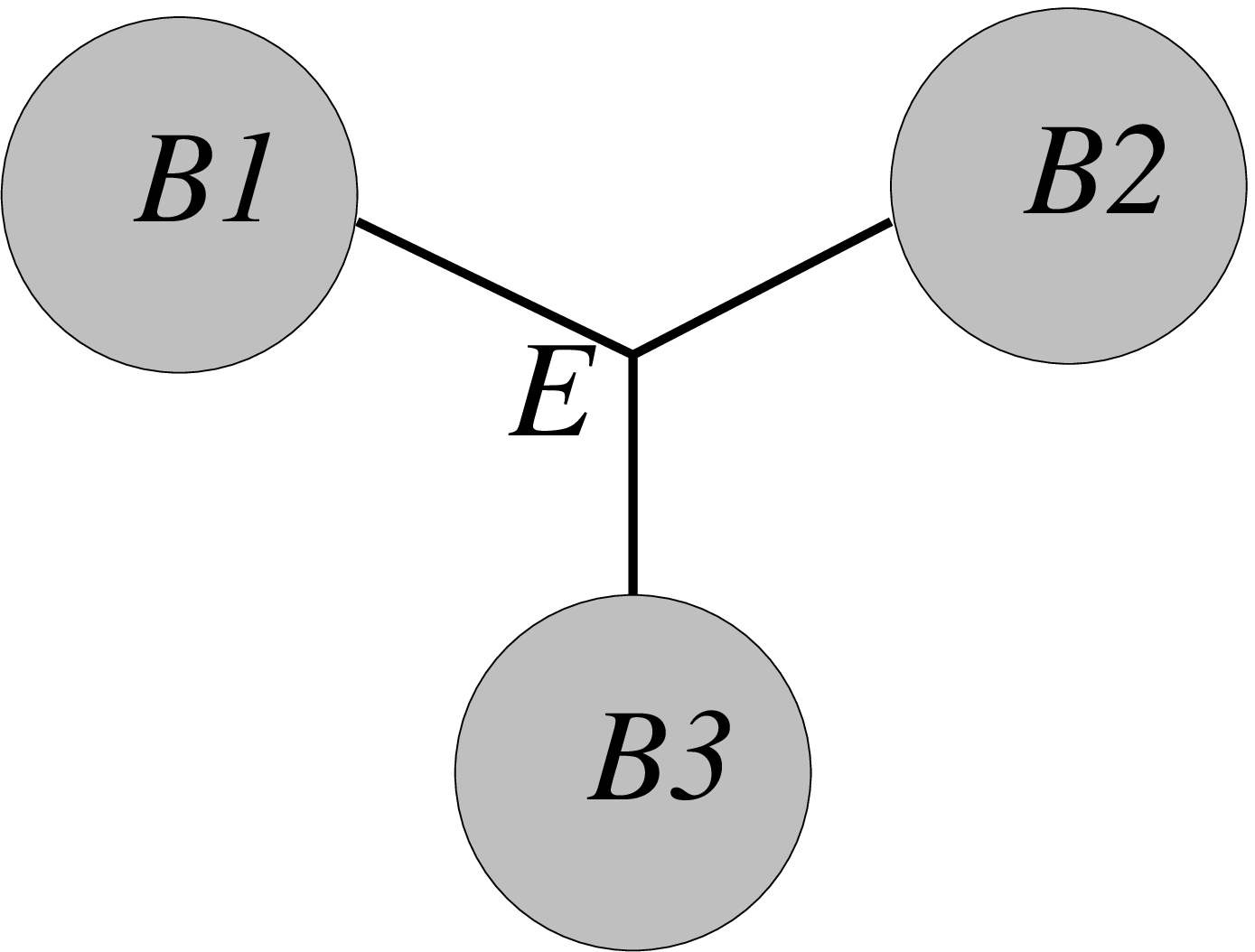}
\end{psfrags}
} \quad\quad \subfigure[sets $\O$ obtained by small perturbation of
$\O^0$.]{\begin{psfrags} \psfrag{O}{$\Omega^n$}
\includegraphics[height=4cm]{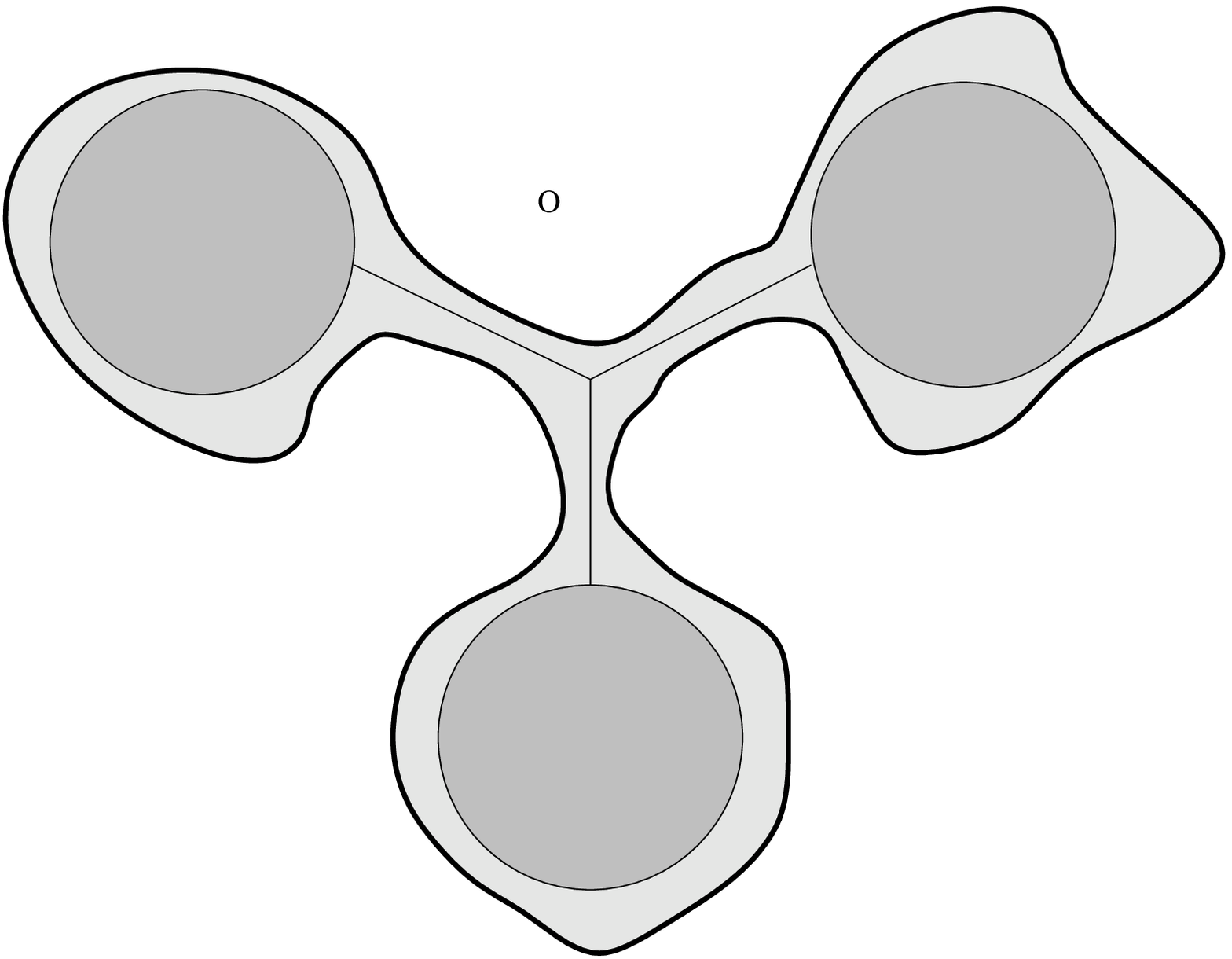}\end{psfrags}}
\caption{}\label{fig:dd}
\end{figure}

In the present paper we show that the geometry of $\Omega$ can
play a crucial role in segregation phenomena.
In line with \cite{efm,mhmm}, where
two populations in planar domains of dumbbell
shape are dealt with, we consider a class of
non convex domains $\Omega^n$, $n\in{\mathbb N}$, essentially composed
by $k$ balls connected by thin corridors, as depicted in Figure
\ref{fig:dd} (see \cite{d} and Section \ref{ipotesi} for the precise definition).
Under the main assumption that the Dirichlet problems on each ball
admit a nondegenerate local minimizer, we are able to prove existence
and uniqueness of positive solutions to the free boundary problem
(\ref{segrego}--\ref{2kvar}), where each component is close to such a
local minimizer (Theorem~1).
In ecological terms this means that if in the unperturbed domain and in absence of
interaction each species lives in a stable configuration, then
strong competition leads to coexistence and segregation of the populations.

Our second result (Theorem 2), concerns the multispecies
Lotka-Volterra system endowed with the Dirichlet null condition. Under
the same topological and nondegeneration assumptions, we first prove
the existence of a positive solution, provided that the competition
parameter $\kappa$ is large enough.
This is obtained by exploiting a degree technique introduced by Dancer \cite{d} to control domain
perturbation in the case of certain nonlinear equations.
Theorem 2 is by itself an interesting
result in the framework of multispecies systems. In fact, in spite of
the rich literature dealing with the case $k=2$ of two populations,
the case of $k\geq 3$ species is less understood. We quote for instance
\cite{eflg, lmn, ln, kl} for three-species competing systems with
cross-diffusion and \cite{ddh,dd} for the Lotka-Volterra model, where
various sufficient conditions for coexistence are provided, depending
on the values of the parameters involved.

Next we perform the asymptotic analysis as $\kappa$ grows to infinity,
and we prove that this solution converges to the unique segregation
state found in Theorem 1.  The biological implication of this result
is now clear: all the species survive under strong competition in a
segregating configuration. Furthermore, as we shall see, they divide
the domain in such a way that the $i$-th species does not invade the
native territory $B_j$ of the other populations. We call this phenomenon
{\it non invading property}.

Both results come from the study of a multispecies system that can be
seen as a generalized Lotka-Volterra model with presence of spatial
barriers localized in the balls.  We shall introduce it in
\eqref{nostro}, after some rigorous definitions and precise statements
of our results.

\subsection{Assumptions and main results}\label{ipotesi}
Let $\Omega^0:=\bigcup_{i=1}^k B_i$ be a finite union of open balls
$\Bi\subset{\bbbr}^N$ such that $\overline{\Bi}$ are disjoint, $i=1,...,k$.
Following \cite{d}, we consider a sequence of domains
$\{\O^n\}_{n\in{\mathbb N}}$ approximating $\O^0$ in the following
sense: there exists a compact zero
measure set $E\subset \R^N$ such that
\begin{align*}
{(\rm i)}\quad &\text{for any compact set } K\subset \Omega^0,\
\Omega^n\supset K\text{ provided $n$ is large};\\
{(\rm ii)}\quad &\text{for any open set } U\supset E\cup\ov{\Omega^0},\
\Omega^n\subset U\text{ provided $n$ is large},
\end{align*}
see Figure \ref{fig:dd}.
Let us fix a bounded smooth domain  $\Omega$ strictly containing
$\ov{\Omega^0}\cup\ov{\O^n}$ for all $n\in{\mathbb N}$.
Notice that if $\tilde\Omega\subset\Omega$ and $u\in H^1_0(\tilde\O)$,
it is possible to extend $u$ to
an element of $H^1_0(\O)$ by defining it to be zero outside of
$\Omega$.
Thus in all the paper we shall think of all our functions as
being in $H^1_0(\O)$.
We will make the
following set of assumptions (for every $i=1,\dots,k$):
\begin{itemize}
\item[\bf(F1)]  $f_i(x,s):\Omega\times \mathbb R\to \mathbb R$ is a
  Carath\'eodory function, it is odd and $C^1$ in
  the variable $s$, uniformly in $x$;
 \smallskip
\item[\bf(F2)]  $|f'_i(x,s)|=O(|s|^{q-1})$ for large $|s|$, uniformly
  in $x$ for some $q<\frac{N+2}{N-2}$ ($q<\infty$ if $N=2$),\medskip
\end{itemize}
where $f'_i(x,s)=\partial_s f_i(x,s)$.
Furthermore,
for any $i=1,\dots,k$, we assume that the problem
\begin{equation}
\label{Ball}
\begin{cases}
 -\D u=f_i(x,u),&\text{ in } B_i,\\
 \quad\;\; u=0,&\text{ on }\partial B_i,
 \end{cases}
\end{equation}
admits a positive solution $u_i^0\in H^1_0(B_i)\cap L^\infty(B_i)$
which is nondegenerate in the following sense:
\begin{itemize}
\item[\bf(ND)]
there exists
$\varepsilon>0$ such that
\begin{equation*}
\int_{B_i} (|\nabla
w|^2-f_i'(x,u_i^0)w^2)dx\geq\varepsilon\int_{B_i}|\nabla w|^2dx,
\end{equation*}
for every $w\in H^1_0(B_i)$ and for every $i$.
\end{itemize}
Note that this implies that the linearized problem at $u_i^0$
\begin{equation*}
-\D v-f_i'(x,u^0_i)v=0,\quad v\in H^1_0(B_i),
\end{equation*}
has only the trivial solution.  This is precisely the assumption used
in \cite{d}.  Condition {\bf(ND)} is stronger and essentially means
that $u_i^0$ is a local nondegenerate minimizer of the energy on $B_i$
(see also \cite{ctv-var}).  As a model for $f_i$ we can consider
logistic type nonlinearities $f_i(x,s)=\lambda(s-|s|^{p-1}s)$, $p>1$.
It is well known that if $\lambda>\lambda_1(B_i)$ (being
$\lambda_1(B_i)$ the first eigenvalue of $-\Delta$ in $B_i$ with
homogeneous Dirichlet boundary conditions) the elliptic problem
\eqref{Ball} has a unique positive solution which is a nondegenerate
global minimum for the energy.

Since $u_i^0$ is an isolated solution to \eqref{Ball}, the parameter
$\delta>0$ appearing in all the paper will be assumed small enough to
ensure that, for all $i$:
$$
\text{\it if $u_i\in H^1_0(B_i)$ is a solution to \eqref{Ball}
  such that $\|u_i-u_i^0\|_{H^1_0(B_i)}\leq \delta$, then $u_i\equiv
  u_i^0$.}
$$
We shall also denote as $U^0$ the $k$-tuple $(u^0_1,\dots,u^0_k)$ and as
$U=(u_1,\dots,u_k)$ generic $k$-tuples in
$\big(H^1_0(\Omega^n)\big)^k$. Let us clarify the meaning of solution to
differential inequalities \eqref{2kvar} in the following definition.
\begin{defin}\label{def:dise}
A solution to \eqref{2kvar} is an $k$-tuple $U=(u_1,\dots,u_k)$ such that,
for every $i=1,\dots,k$ and $\phi\in H^1_0(\Omega)$,
$\phi\geq 0$ a.e. in $\Omega$, there holds
$$
\int_\Omega \nabla u_i(x)\cdot \nabla\phi(x)\,dx\leq \int_\Omega
f_i(x,u_i(x))\phi(x)\,dx
$$
and
$$
\int_\Omega \nabla \widehat u_i(x)\cdot \nabla\phi(x)\,dx\geq
\int_\Omega \widehat f_i(x,\widehat u_i(x))\phi(x)\,dx.
$$
\end{defin}

\noindent The following theorem ensures the existence of a unique
segregated solution to \eqref{2kvar} in the perturbed domain
$\Omega^n$ which is $H^1$-close to $U^0$.

\medskip
\noindent
\textbf{Theorem 1.}
\textit{Let us define
\begin{equation*}
\mathcal S(\Omega^n)=\left\{
  \begin{array}{l}
    (u_1,\dots, u_k)\in \big(H^1_0(\Omega^n)\big)^k\;:
    \;u_i\geq 0,\; u_i\cdot u_j=0,\text{ if }i\neq j\\
    -\Delta u_i\leq f_i(x,u_i),\; -\Delta \widehat u_i
    \geq \widehat f_i(x,\widehat u_i),\text{ in } \Omega^n,\;
    i=1,\dots,k
  \end{array}\right\}.
\end{equation*}
Then, there exists $\delta>0$ such that, for any $n$ sufficiently
large, the class $\mathcal S(\Omega^n)$ contains an element
$U=(u_1,\dots, u_k)\in \big(H^1_0(\Omega^n)\big)^k$ such that
$\|u_i-u_i^0\|_{H_0^1(\Omega^n)}<\delta$.
 Moreover,  $u_i\equiv 0$ in $B_j$
for all $j\neq i$ and $U$ is the unique element of
$S(\Omega^n)$ such that $\|U-U^0\|_{(H^1_0(\Omega^n))^k}<\delta$,
where $U^0$.}

\smallskip
As we already announced, the proof of Theorem 1 relies on a careful analysis
of the following  auxiliary system:
\begin{equation}
\label{nostro}
\begin{cases}
 -\D u_i=f_i(x,u_i)\dis- \varkappa u_i\sum_{j\neq i}u_j\dis- \varkappa
 u_i\sum_{j\neq i}u_j^0
 \dis- \varkappa u^0_i\sum_{j\neq i}u_j,
  &\text{ in }\Omega^n,\\
 \quad\;\; u_i=0,&\text{ on }\partial\Omega^n,
 \end{cases}
\end{equation}
for $i=1,\dots,k$.

This system can be seen as a modification of the Lotka-Volterra model,
through linear terms which are localized in the single balls: this
feature will be crucial in order to obtain solutions with the {\it
  non-invading} property.  Notice that, due to the presence of the
barriers, systems \eqref{nostro} lack of the maximum principle, so
that we cannot ensure the positivity of its solutions nor even, by
now, the competitive character of the model.  As we shall see, this
will cause some technical difficulties.

Nonetheless, by careful energy estimates and eigenvalue theory, the
system at fixed $\kappa$ will be shown to be suitably nondegenerate on
$\Omega^n$, if $n$ is large enough. This will allow the application of
the degree technique introduced in \cite{d} to control domain
perturbation, and as a result
we will obtain the existence of a solution $U^\kappa$ of the system,
which is close to $U^0$.  The major feature of this approach is that
the whole procedure turns out to be {\it uniform} with respect to
$\vk$. This uniformity will allow to perform successfully the
asymptotic analysis of the solutions to the auxiliary system as the
competition parameter goes to infinity.  The final result can be
collected in the following form.

\medskip
\noindent \textbf{Theorem 2.}
\textit{There exists $\delta>0$ such that, for any $\vk$ and $n$
  sufficiently large, both the Lotka-Volterra system \eqref{modelLK}
  with Dirichlet boundary conditions \eqref{dirich} on $\Omega^n$ and
  the modified model \eqref{nostro}, admit a solution
  $U^{\vk}=(u_1^{\vk},\dots,u_k^{\vk})\in \big(H^1_0(\Omega^n)\big)^k$
  such that $\|U^{\vk}- U^0\|_{(H^1_0(\Omega^n))^k}<\delta$, and, in
  the case of \eqref{modelLK}, $U^{\vk}$~is strictly positive.  Furthermore, as
  $\vk\to\infty$, $u_i^{\varkappa}\to u_i$ strongly in
  $H^1(\Omega^n)$, where the $k$-tuple $U=(u_1,\dots,u_k)$ is the
  unique element in $\mathcal S(\Omega^n)$ close to~$U^0$.}

\subsection{Plan of the paper.}
In Section \ref{preliminar} we establish some preliminary facts that
shall be used throughout the paper, in particular we discuss the
nondegeneracy of the problems in $\Omega^0$.  Section \ref{asintotic}
is devoted to perform the asymptotic analysis of the solutions to the
auxiliary system as $\vk\to\infty$.  In Section \ref{unicos} we prove
the uniqueness of the solution to \eqref{2kvar} close to $U^0$, as
stated in the uniqueness part of Theorem 1.  Section \ref{perturb} is
devoted to the proof of the existence of a solution $U$ close to $U^0$
for system \eqref{nostro}, when the domain is close enough to
$\Omega^0$ and the competition is large.  We conclude section
\ref{perturb} by presenting the proofs of Theorems 1 and 2 and giving
some final remarks. A final appendix collects some technical proofs
and lemmas used throughout the paper.

\section{Preliminary results}\label{preliminar}
In this section we further modify the Lotka-Volterra system in order
to ensure sign conditions and boundedness of its solutions.
Furthermore we derive by condition {\bf (ND)} the main nondegeneracy
properties holding in the unperturbed set $\Omega^0$.

Let us consider the following system:
\begin{equation}
\label{problpiu}
\begin{cases}
 -\D u_i=f_i(x,[u_i+u_i^0]^+-u_i^0)\dis-
 \varkappa [u_i+u_i^0]^+\sum_{j\neq i}[u_j+u_j^0]^+,
  &\text{ in }\tilde\Omega,\\
 \quad\;\; u_i=0,&\text{ on }\partial \tilde\Omega,
 \end{cases}
\end{equation}
for $i=1,\dots,k$, where
$\Omega^0\subset\tilde\O\subseteq\Omega$.

Here and throughout the paper the symbol $[t(x)]^+$ will denote the
positive part of $t$, namely $[t(x)]^+=\max\{t(x),0\}$.  The
motivation for this choice is contained in the following lemma
\begin{lemma}\label{l:mp}
Let $U=(u_1,\dots,u_k)$ be a solution of system
\eqref{problpiu}.
Then
for all $i=1,\dots,k$, $u_i\geq -u_i^0$ in  $\tilde\Omega$. In particular
they solve \eqref{nostro} and
satisfy $u_i(x)\geq 0$ for $x\in B_j$ when $j\neq i$.
\end{lemma}
\proof
Let $v_i=u_i+u_i^0$ with its equation
$$-\D v_i=f_i(x,[u_i+u_i^0]^+-u_i^0)+f_i(x,u_i^0)\dis-
\varkappa [v_i]^+\sum_{j\neq i}[v_j]^+.$$ Then it suffices to test the
equation by $-[u_i+u_i^0]^-$, recalling that $f_i$ is odd.
\endproof
\begin{rem}\label{LMpiu}{\rm Notice that the original system
    \eqref{modelLK} can be recovered in this model by the formal
    identification $u_i^0\equiv 0$.  In particular by Lemma \ref{l:mp}
    it turns out that the solutions obtained with $f_i(x,[u_i]^+)$
    instead of $f_i(x,u_i)$ satisfy $u_i\geq 0$ for all $i$, and thus
    are nonnegative solutions for the Lotka-Volterra system
    \eqref{modelLK}.}
\end{rem}
\subsection{Differential inequalities.}\label{prelim}
Let $(u_1,\dots,u_k)$ be a solution to \eqref{problpiu}. Since by
Lemma \ref{l:mp} each $u_i$ satisfies $u_i+u_i^0\geq 0$, the coupling
term has negative sign and we immediately have \neweq{sotto}-\Delta
u_i\leq f_i(x,u_i).\endeq Furthermore, by a straightforward
calculation we obtain an opposite differential inequality for
$\hatu_i$: \neweq{sopra}-\Delta \hatu_i\geq
f_i(x,u_i)-\displaystyle\sum_{j\neq i} f_j(x,u_j)=\widehat
f_i(x,\widehat u_i).\endeq It turns out that the solutions of
\eqref{problpiu} satisfy the differential inequalities \eqref{2kvar}.

\subsection{Uniform $L^\infty$ bounds.}\label{sec:unif-linfty-bounds}
We now suitably modify $f_i$ in order to ensure that the solutions of
the new system are bounded in $L^\infty$. This is based on the
following result due to Dancer \cite{d}: \textit{for $n$ sufficiently
  large, the problem
\begin{equation*}
\begin{cases}
 -\D u=f_i(x,u),&\text{ in } \O^n,\\
 \quad\;\; u=0,&\text{ on }\partial \O^n,
 \end{cases}
\end{equation*}
admits a positive solution $\phi^n_i\in H^1_0(\O^n)$ which is close to
$\sum u_i^0$ in some $L^r(\O)$ $(r>1)$.}
Let $n_0$ large such that $\O^n\subset \O^{n_0}$ for all $n\geq n_0$
and denote $\phi_i^{n_0}$ simply $\phi_i$.
Let us define
\begin{equation*}\tilde f_i(x,s)=
\begin{cases}
 f_i(x,s),&\text{ if } s\leq \phi_i(x),\\
 f_i(x,\phi_i(x)),&\text{ if } s>\phi_i(x).
 \end{cases}
\end{equation*}
\begin{lemma}\label{l:ublinf}
  For $n\geq n_0$, let $u_i\in H^1_0(\O^n)$ such that $-\Delta u_i\leq
  \tilde f_i(x,u_i)$ in $\O^n$. Then $u_i\leq \phi_i$ a.e. in $\O$.
\end{lemma}
\proof
Summing up the differential inequalities for $\phi_i$ and $u_i$ it holds
\begin{equation*}
\begin{cases}
 -\D (\phi_i-u_i)\geq f_i(x,\phi_i)-\tilde f_i(x,u_i),&\text{ in } \O^n,\\
 \quad\;\; \phi_i-u_i\geq 0,&\text{ on }\partial \O^n.
 \end{cases}
\end{equation*}
Set $\omega=\{x\in\O^n: \phi_i< u_i\}$: note that $\omega$ is strictly
contained in $\O^n$ by the boundary conditions. Hence by  testing the
first inequality  with $-(\phi_i-u_i)^-$ we obtain
$$
\int_\omega |\nabla (\phi_i-u_i)^-|^2\,dx\leq
-\int_\omega\big(f_i(x,\phi_i)-\tilde f_i(x,u_i)\big)(\phi_i-u_i)^-\,dx=0,
$$
which implies $(\phi_i-u_i)^-=0$ and  so $\phi_i\geq u_i$.
\endproof

\noindent
As a consequence, any solution $(u_1,\dots,u_k)$ either to \eqref{problpiu}
or to \eqref{2kvar}, with $\tilde f_i$ instead of $f_i$ in $\O^n$,
satisfies
\begin{equation}\label{pinzo}
-u_i^0\leq u_i\leq
\phi_i,\quad i=1,\dots,k.
\end{equation}
In particular, any solution with $u_i$ close to
$u_i^0$
will be a true solution of the original problem with $f_i$.
Moreover, by virtue of a classical strong maximum principle and Harnack's
inequality, any solution $(u_1,\dots,u_k)$ to \eqref{problpiu} satisfies
either $u_i\equiv \phi_i$, or $u_i\equiv -u_i^0$, or
$-u_i^0<u_i<\phi_i$.

\medskip\noindent
{\bf Notations.} Throughout all the paper we shall work with $\tilde
f_i$ instead of $f_i$, denoting $\tilde f_i$ simply by $f_i$.
$\supp(u_i)$ will denote the set $\{u_i>0\}$.

\subsection{Nondegeneracy in $\Omega^0$.}
Let us now consider \eqref{problpiu} in $\Omega^0$: then we
immediately realize that $U^0=(u_1^0,\dots,u_k^0)$ is a solution of
the problem.  Furthermore, it comes from {\bf (ND)} that $U^0$ is an
isolated solution, uniformly in $\vk$.
\begin{teo}\label{c:isol}
  There exists $\bar\vk>0$ and $\delta>0$ such that if $U^\varkappa$
  is a solution of \eqref{problpiu} in $\Omega^0$ such that
  $\|U^\varkappa-U^0\|<\delta$ in $\big(H^1_0(\UB)\big)^k$, then
  $U^\varkappa\equiv U^0$ for all $\vk\geq\bar\vk$.
\end{teo}
An analogous result holds for the solutions to
\eqref{segrego}-\eqref{2kvar}, thanks to the following sign condition
prescribed by the validity of \eqref{sopra}.
\begin{lemma}\label{sopravuota}
  Let $(u_1,\dots, u_k)$ be solution of $-\Delta \widehat u_i\geq
  \widehat f(x,\widehat u_i)$ in $B_i$ for some $i$.  Assume that
  $u_i\cdot u_j=0$ if $i\neq j$ and that
  $\|u_j-u_j^0\|_{H^1_0(\Omega)}\leq \delta$ for all $j=1,\dots,k$.
  Then, if $\delta$ is small enough, $\widehat u_i\geq0$.
\end{lemma}

\begin{teo}\label{essevuota}
  Let $(u_1,\dots, u_k)\in \mathcal S(\Omega^0)$ such that
  $\|u_i-u_i^0\|_{H^1_0(\Omega^0)}\leq \delta$ for all $i=1,\dots, k$.
  Then, if $\delta$ is small enough, $u_i\equiv u_i^0$ for all
  $i=1,\dots,k$.
\end{teo}
All these results are crucial in what follows, but since the proofs
are somewhat technical, we postpone them in the Appendix.

\section{Asymptotic analysis as $\vk\to\infty$}\label{asintotic}
This section is devoted to establish the link between the population
systems and the original set of differential inequalities
\eqref{2kvar}. To this aim, throughout the whole section let
$\delta>0$ and assume that {\it there exists $(u_1^\vk,\dots,
  u_k^\vk)$ solution to \eqref{problpiu} such that $\|
  u_i^\vk-u_i^0\|_{H^1(\Omega)}\leq \delta$ for all large $\vk$.}
Our main result is
\begin{teo}\label{limite}
  Let $\tilde\O$ be a connected domain such that
  $\Omega^0\subset\tilde\O\subseteq\Omega$.  For each $\varkappa$ let
  $U^\varkappa=(u_1^\varkappa,...,u_k^{\varkappa})$ be a solution of
  \eqref{problpiu} in $\tilde\O$ such that
  $\|U^{\vk}-U^0\|_{(H^1_0(\O))^k}<\delta.$ Then, if $\delta$ is small
  enough, there exists $U\in\big(H^1_0(\tilde\O)\big)^k$ such that,
  for all $i=1,\dots,k$:
  \begin{itemize}
  \item[(i)] up to subsequences, $u_i^{\varkappa}\to u_i$
    strongly in $H^1$ as $\vk\to\infty$,
  \item[(ii)] $u_i>0$ in $B_i$,
  \item[(iii)] if $i\neq j$ then $u_i=0$ a.e. in $B_j$,
  \item[(iv)] if $i\neq j$ then $u_i\cdot u_j=0$ a.e. in $\O$,
  \item[(v)] $(u_1,\dots,u_k)$ satisfy the differential inequalities
    \eqref{2kvar}.
    \end{itemize}
\end{teo}
The proof of this fact is obtained through the next Lemmas
\ref{limite}, \ref{positiva} and \ref{limiteforte}.
\begin{lemma}\label{limiteweak}
  Under the assumptions of Theorem \ref{limite}, if $\delta$ is small
  enough there exists $U\in\big(H^1_0(\tilde\O)\big)^k$ such that, for
  all $i=1,\dots,k$:
  \begin{itemize}
  \item[(i)] up to subsequences, $u_i^{\varkappa}\rightharpoonup u_i$
    weakly in $H^1(\O)$ as $\vk\to\infty$,
  \item[(ii)] $u_i\geq 0$ in $\Omega\setminus B_i$,
  \item[(iii)] if $i\neq j$ then $u_i=0$ a.e. in $B_j$,
  \item[(iv)] if $i\neq j$ then $u_i\cdot u_j=0$ a.e. in $\O$,
  \item[(v)] $(u_1,\dots,u_k)$ satisfy the differential inequalities
    \eqref{2kvar}.
      \end{itemize}
\end{lemma}
\proof Since $U^\varkappa$ is bounded in $(H^1(\O))^k$ by assumption,
we immediately obtain the existence of a weak limit $U$ such that, up
to subsequences, $u_i^\varkappa\rightharpoonup u_i$ in $H^1(\O)$.
Since each $u_i^\varkappa$ is positive on $B_j$ when $j\neq i$ by
Lemma \ref{l:mp}, property $(ii)$ comes from almost everywhere
pointwise convergence. Furthermore, the differential inequalities
\eqref{sotto} and \eqref{sopra} for $u_i^\varkappa$ pass to the weak
limit, so $(v)$ is already proved.  Let us discuss properties (iii)
and (iv). By testing \eqref{nostro} times $u_i^\vk+u_i^0$ we obtain
\begin{eqnarray*}
  \vk\int_{\Omega}(u^\vk_i+u_i^0)^2\sum_{j\neq
    i}(u^\vk_j+u_j^0)\quad\text{is bounded uniformly in }\vk,
\end{eqnarray*}
hence, since $u^\vk_j+u_j^0\geq0$ for all $j$,
$$
\int_{\Omega}(u^\vk_i+u_i^0)^2\sum_{j\neq i}(u^\vk_j+u_j^0)\to
0,\qquad \text{as }\vk\to\infty.
$$
Passing to the limit for $U^\vk\rightharpoonup U$ we obtain, for all
$i\neq j$, $i,j=1,...,k$
\begin{equation}
\label{disgiunti}
u_i(x)\cdot u_j(x)+u_i^0(x)\cdot u_j(x) +u_i(x)\cdot
u_j^0(x)=0,\qquad\forall x\in\Omega.
\end{equation}
Let $x\in\tilde\Omega\setminus \cup B_i$: then \eqref{disgiunti}
ensures $u_i(x)\cdot u_j(x)=0$ for all $i\neq j$.

{\it Claim. If $x\in B_i$  then $u_j(x)=0$ for all $j\neq i$.}

Let $x\in B_i$ for some fixed $i$. If $u_i(x)=0$ then
\eqref{disgiunti} becomes $u_i^0(x)\cdot u_j(x)=0$ and hence
$u_j(x)=0$, for all $j\neq i$. If $u_i(x)>-u_i^0(x)$ we have
$u_j(x)\big(u_i(x)+u_i^0(x)\big)=0$ implying again $u_j(x)=0$ for all
$j\neq i$.  Finally, let $u_i(x)=-u_i^0(x)$. Since $u_j\cdot u_h=0$ in
$B_i$ for all $j\neq h$, $j,h\neq i$, then there exists at most one
index different from $i$ (say $j$) where $u_j(x)>0$. Set $\omega_j$
the connected component of $\{u_j>0\}$ which is contained in the set
$\{y\in B_i:\, u_i(y)=-u_i^0(y)\}$, and such that $x\in \omega_j$.
Then by \eqref{sopra}, since $u_h=0$ in $\omega_j$ for all $h\neq
i,j$, we have
$$
-\Delta (u_i-u_j)\geq f_i(\cdot,u_i)-f_j(\cdot,u_j)
\qquad\text{ in }\omega_j.
$$
Adding $-\Delta u_i^0=f_i(\cdot,u_i^0)$ we get
$$
-\Delta (u_i-u_j+u_i^0)\geq f_i(\cdot,u_i)-f_j(\cdot,u_j)
+f_i(\cdot,u_i^0),\qquad\text{ in }\omega_j.
$$
Test this equation times $-[u_i-u_j+u_i^0]^-.$ Note that
$[u_i-u_j+u_i^0]^-\equiv u_j|_{\omega_j}$, providing
$$
\int_{\omega_j}|\nabla u_j|^2\leq \int_{\omega_j}f_j(x,u_j)u_j\leq
\|f_j(\cdot,u_j)\|_{L^{N/2}(\omega_j)}\|u_j\|^2_{L^{2^*}(\omega_j)}.
$$
By {\bf (F1)}, {\bf (F2)}, since $\|u_j\|_{H^1(\omega_j)}\leq \delta$,
we have $\|f_j(\cdot,u_j)\|_{L^{N/2}(\omega_j)}\leq C\,\delta$ which
implies $u_j\equiv 0$ in $\omega_j$ if $\delta$ is small enough,
giving rise to a contradiction. This proves the claim.
\endproof
\begin{lemma}\label{positiva}
  Under the assumptions of Theorem \ref{limite}, if $\delta$ is small
  enough, then $u_i>0$ in the whole of $B_i$. In particular
  $B_i\subset \supp(u_i)$.
\end{lemma}
\proof By Theorem \ref{limiteweak}, we already know that $u_i\geq 0$
in $\O\setminus B_i$. Furthermore, $u_i\cdot u_j=0$ for $i\neq j$ in
$\Omega$ and $(u_1,\dots,u_k)$ satisfies \eqref{sopra}.  Therefore
Lemma \ref{sopravuota} yields $\widehat u_i=u_i-\sum_{j\neq i}u_j\geq
0$ in $B_i$.  Since by (ii) of Theorem \ref{limite} $u_j\geq 0$ in
$B_i$, then $u_i\geq 0$ in $B_i$.  Hence $u_i\geq 0$ in $\Omega$ and
it is not identically null by its closeness to $u_i^0$; the strict
positivity now comes from the Harnack inequality.
\endproof
\begin{lemma}\label{limiteforte}
  Under the assumption of Theorem \ref{limite}, the convergence
  $u_i^{\varkappa}\to u_i$ is strong in~$H^1_0(\O)$ (up to
  subsequences), where $U=(u_1,\dots,u_k)$ is as in Lemma
  \ref{limiteweak}.
\end{lemma}
\proof
In order to prove the strong convergence of $u_i^{\vk}$ to $u_i$ in
$H^1_0(\Omega)$, let us consider the functions $\hat
u_i=u_i-\sum_{j\neq i} u_j$, which satisfy the inequality \eqref{sopra}
in $\tilde\Omega$.
Since from Lemma \ref{positiva} $u_i\geq 0$, testing \eqref{sopra}
 with $u_i$ we obtain
\begin{equation}\label{eq:vit1}
\int_{\Omega}u_i f_i(x,u_i) \leq \int_{\Omega}|\nabla u_i|^2.
\end{equation}
Testing
\begin{equation*}
-\Delta  u_i^{\vk}\leq f_i(x,u_i^{\vk}),\quad\text{in }\tilde\Omega,
\end{equation*}
with $u_i^{\vk}+u_i^0$ (which is positive in view of Lemma \ref{l:mp})
we have
\begin{equation}\label{eq:vit}
\int_{\Omega}\nabla
u_i^{\vk}\cdot\nabla u_i^0+\int_{\Omega}|\nabla u_i^{\vk}|^2
\leq \int_{\Omega}u_i^0
f_i(x,u_i^{\vk})+\int_{\Omega}u_i^{\vk} f_i(x,u_i^{\vk}).
\end{equation}
The uniform $L^{\infty}$-bound provided in Section
\ref{sec:unif-linfty-bounds} and the Dominated Convergence Theorem
allows to pass to the limit in \eqref{eq:vit}  thus obtaining
\begin{equation*}
\int_{B_i}\nabla u_i^0\cdot \nabla u_i+\limsup_{\vk\to\infty}
\int_{\Omega}|\nabla u_i^{\vk}|^2\leq \int_{B_i}u_i^0
f_i(\cdot,u_i)+\int_{\Omega}u_i f_i(\cdot,u_i).
\end{equation*}
Since by Theorem \ref{limite},  $u_i$ solves $-\Delta u_i=f_i(x,u_i)$ in
$B_i$, testing with $u_i^0$ we have $\int_{B_i}\nabla
u_i\cdot\nabla u_i^0= \int_{B_i}u_i^0f_i(x,u_i)$, thus implying
\begin{equation}\label{eq:vit3}
\limsup_{\vk\to\infty}
\int_{\Omega}|\nabla u_i^{\vk}|^2\leq \int_{\Omega}u_i f_i(\cdot,u_i).
\end{equation}
Now \eqref{eq:vit1}, \eqref{eq:vit3}, and the lower semi-continuity of
the norms yield
$$
\lim_{\vk\to\infty}
\int_{\Omega}|\nabla u_i^{\vk}|^2=\int_{\Omega}|\nabla u_i|^2.
$$
The strong convergence follows easily from weak convergence and
convergence of norms.~\endproof
\begin{rem}\label{LMasynt}{\rm Notice that the above analysis can be
    performed also for the Lotka-Volterra system \eqref{modelLK}, with
    some differences. In particular, following the proof of Lemma
    \ref{limite}, the segregation property $(iv)$ immediately follows
    by \eqref{disgiunti} which reduces in this case to $u_i\cdot
    u_j=0$.  On the contrary we cannot prove, at the moment, the
    noninvading property $(iii)$.}
\end{rem}

\section{Uniqueness of the asymptotic limit.}\label{unicos}
As in the previous section, let us here assume that the system
\eqref{problpiu} does have a solution on $\Omega^n$ for all $\vk$
large. Our goal now consists in proving that the class $\mathcal
S(\Omega^n)$ contains one single element which is close to $U^0$;
it is worth noticing that $U^0$ does not belong to $\mathcal
S(\Omega^n)$, since the differential inequalities involving the hat
operation cannot hold outside $\Omega^0$.
\begin{teo}\label{uniqueness}
  For $\delta$  sufficiently small and  $n$ sufficiently large,
the class $\mathcal
  S(\Omega^n)$ has at most a unique element $U$ such that
  $\|U-U^0\|_{((H^1_0(\Omega))^k}<\delta $.
\end{teo}
\proof
 By Theorem \ref{limite}, let $U^n\in \mathcal S(\Omega^n)$ be the
 asymptotic limit of the solutions to \eqref{problpiu}, so that $U^n$
 enjoys the noninvading property. Now
assume by contradiction the existence of  $V^n\in\mathcal
S(\Omega^n)$, such that  $U^n\neq V^n$.

{\it Claim 1.
Letting $n\to\infty$, both $U^n\to U^0$ and $V^n\to U^0$ weakly in
$H^1_0(\Omega)$ (hence strongly in $L^p(\O)$ for
all $1\leq p< 2^*$).}

It suffices to prove the claim for $U^n$. Since $U^n$ is bounded in
$(H^1(\O))^k$, there exists $U\in (H^1(\Omega))^k$ such that
$u_i^n\rightharpoonup u_i$ weakly in $H^1(\O)$, strongly in all
$L^p(\Omega)$ with subcritical $p$.  We are going to prove that $U\in
\mathcal S(\O^0)$ so that $U\equiv U^0$ in light of Lemma
\ref{essevuota}.  To this aim notice that the differential
inequalities characterizing $\mathcal S(\O^0)$ are satisfied by $u_i^n$
for all $n$, hence they pass to the weak limit.  It remains to prove
that $u_i\in H^1_0(\O^0)$. To this aim notice that, for all open sets
${\mathcal V}$ containing $\ov\O^0\cup E$, we have that
$$
\mathop{\rm supp}u_i^n\subset \O^{n}\subset{\mathcal V},
$$
provided that $n$ is sufficiently large.
Hence
$$
\mathop{\rm supp}u_i\subset {\mathcal V}\quad\text{for all  open sets
  ${\mathcal V}$ containing $\ov\O^0\cup E$},
$$
which implies that $u_i=0$ a.e. in $\O\setminus\big(\ov\O^0\cup
E\big)$.  Since $\partial\Omega\cup E$ has measure zero, $u_i=0$ on
$\O\setminus\ov\O^0$, and the smoothness of $\partial \O^0$ ensures that
$u_i\in H^1_0(\O^0)$ (see \cite{GilTr}).

Let us now start the argument that will lead to a contradiction.
By setting $\omega_i^n=\{u_i^n>0\}$, we have
$\widehat u_i^n=u_i^n$ in $\omega_i^n$ and  the following hold:
\begin{align*}
-\Delta  u_i^n=f_i(x, u_i^n),\hspace{1cm}&\text{in }\omega_i^n,\\
-\Delta  v_i^n\leq f_i(x, v_i^n),\hspace{1cm}&\text{in }\omega_i^n.
\end{align*}
If we now consider
$$
w_i^n=\frac{ v_i^n-u_i^n}{\| V^n-U^n\|_{L^2(\Omega)}},
$$
we have
\begin{equation}
  \label{uniq1}
  \left\{\begin{array}{l}
      -\Delta w_i^n \leq a_i^n(x)w_i^n,\hspace{1cm}\text{in }\omega_i^n,\\
      -\Delta \widehat w_i^n\geq b_i^n(x)\widehat w_i^n,
      \hspace{1cm}\text{in }\omega_i^n,
    \end{array}
\right.
\end{equation}
where $a_i^n(x)=\frac{f_i(x,v_i^n)-f_i(x,u_i^n)}{v_i^n-u_i^n}$ and
$b_i^n(x)=\frac{\widehat f_i(x,\widehat v_i^n)-f_i(x,u_i^n)}{\widehat
  v_i^n-u_i^n}$.
Notice that $a_i^n(x)\in L^\infty$ independently of $n$ in light of
the a priori estimates in Remark \ref{pinzo} and Lemma \ref{l:ublinf}
and since $f'_i(\cdot,0)$ is bounded. We assert that this is true for
the second quotient too.  To see this, remember that $v_i^n\cdot
v_j^n=0$ in $\Omega$ and notice that
$b_i^n(x)=\frac{f_i(x,v_i^n(x))-f_i(x,u_i^n(x))}{v_i^n(x)-u_i^n(x)}$
for $x\in\Omega$ such that $v_i^n(x)>0$. On the other side, if
$v_j(x)>0$ for some $j\neq i$, then
$b_i^n(x)=\frac{f_i(x,u_i^n(x))+f_j(x,v_j^n(x))}{u_i^n(x)+v_j^n(x)}$.
Hence the same argument used to estimate $a_i^n(x)$ provides an
$L^\infty$ control for $b_i^n$, uniformly in $n$.  As a consequence,
by testing the differential inequalities in \eqref{uniq1} with
$[w_i^n]^+$ and $-[\widehat w_i^n]^-$ respectively, we easily obtain
that $w_i^n$ is bounded in $H^1(\Omega)$. Since this is true for all
$i=1,\dots,k$, there exists $W=(w_1,\dots,w_k)\in (H^1(\Omega))^k$
such that $w_i^n\rightharpoonup w_i$ weakly in $H^1(\Omega)$, strongly
in $L^2$ so that $w_i\neq 0$ for some $i$.

{\it Claim 2. $w_i\in H^1_0(B_i)$ and $-\Delta w_i\leq
  f_i'(x,u_i^0)w_i$ in $B_i$.}

Reasoning as in Claim 1, we can easily prove  that $w_i\in H^1_0(B_i)$.

Let $\phi\geq0$ such that $\phi\in C^\infty_0(B_i)$.
Since by Theorem \ref{positiva} we know that $B_i\subset \omega_i^n$, we can
test the first inequality in \eqref{uniq1} with $\phi$, providing
$$\int_{B_i}\nabla w^n_i\nabla\phi-a_i^n(x) w^n_i\phi\leq 0.$$
By the strong convergence of $U^n$ and $V^n$ to $U^0$ in $L^p(\O)$ for
all $1\leq p< 2^*$ and by the
 continuity of the Nemytskij operator
$f_i':L^{ N(q-1)/2}(\O)\to L^{N/2}(\O)$, see \eqref{eq:9},  it is easy to
realize that  $a_i^n\to f_i'(\cdot,u_i^0)$ in $L^{N/2}(\Omega)$  as
$n\to\infty$. Hence we can pass to the limit and we find
\begin{equation}\label{eq:11}
\int_{B_i}\nabla w_i\nabla\phi-f_i'(x,u_i^0)w_i\phi\leq 0.
\end{equation}

By exploiting the same argument starting by the second inequality in
\eqref{uniq1}, we can prove the opposite inequality, namely $-\Delta
w_i\geq f_i'(x,u_i^0)w_i$ in $B_i$.  To this aim, we first notice
that, setting $A_n^{i,j}:=\{x\in B_i: v_j^n(x)>0\}$, there holds
$\lim_{n\to+\infty}\mu(A_n^j)=0$ if $j\neq i$. Indeed, since
$A_n^{i,j}\subset\{x\in\Omega\setminus B_j: v_i^n(x)=0\}$ for every
$j\neq i$, we have that
\begin{align*}
o(1)&=\int_{B_i}|u_i^0-v_i^n|^2\,dx\geq
\int_{\{x\in B_i:v_i^n(x)=0\}}|u_i^0|^2\,dx\geq \int_{A_n^{i,j}}|u_i^0|^2\,dx
\end{align*}
as $n\to+\infty$. Due to the absolute continuity of the Lebesgue
measure $\mu$ with respect to the measure $A\mapsto\int_A
|u_i^0|^2\,dx$ in $B_i$, we deduce that
\begin{equation}\label{eq:12}
\lim_{n\to+\infty}\mu (A_n^{i,j})=0.
\end{equation}
In $B_i$ we can write $b_i^n$ as
\begin{align*}
b_i^n(x)=&
\frac{f_i(x,v_i^n(x))-f_i(x,u_i^n(x))}{v_i^n(x)-u_i^n(x)
}\chi_{B_i\cap \mathop{\rm supp}v_i^n }(x)\\
&\quad+
\sum_{j\neq i}
\frac{f_i(x,u_i^n(x))+f_j(x,v_j^n(x))}{u_i^n(x)+v_j^n(x)}
\chi_{A_n^{i,j} }(x).
\end{align*}
From $\lim_{n\to+\infty}\mu (A_n^{i,j})=0$, the a priori estimates in
Remark \ref{pinzo} and Lemma \ref{l:ublinf}, and since $f'_i(\cdot,0)$
is bounded, we deduce that $b_i^n(x)\to f_i'(x,u_i^0(x))$ for a.e.
$x\in B_i$. From the uniform $L^{\infty}$-boundedness of $b_i^n$ and
the Dominated Convergence Theorem, we conclude that $b_i^n\to
f_i'(\cdot,u_i^0)$ in $L^{N/2}(\Omega)$ as $n\to\infty$.  Hence
testing the second inequality in \eqref{uniq1} with $\phi\in C^\infty_0(B_i)$,
$\phi\geq0$, and passing to the limit as $n\to+\infty$ we obtain that
\begin{equation}\label{eq:13}
\int_{B_i}\nabla \widehat w_i\nabla\phi-f_i'(x,u_i^0)\widehat w_i\phi\geq 0.
\end{equation}
From (\ref{eq:12}) and $L^2$-convergence of $w_j^n$ to $w_j$, for
all $j\neq i$ there holds
$$
\int_{B_i}|w_j|^2\,dx=\lim_{n\to+\infty}\int_{B_i}|w_j^n|^2\,dx=
\lim_{n\to+\infty}\int_{A_n^{i,j}}|w_j^n|^2\,dx=0.
$$
Therefore $w_j=0$ a.e. in $B_i$ for every $j\neq i$, and
\begin{equation}\label{eq:14}
\widehat w_i=w_i\quad\text{in }B_i.
\end{equation}
From (\ref{eq:11}), (\ref{eq:13}), and (\ref{eq:14}), we conclude that
 $w_i$ is a nontrivial solution to the linearized equation
$-\Delta w_i= f_i'(x,u_i^0)w_i$ in $B_i$
with boundary condition $w_i=0$ on $\partial B_i$.
This provides a contradiction with the nondegeneracy assumption {\bf(ND)}.
\endproof

\begin{rem}\label{rem:strong}\rm
We note that the weak $H^1(\O)$-convergence stated in Claim 1 of the
above proof, is actually strong. Indeed from Theorem \ref{limite}(v),
we have that
$$
\|u_i^n\|^2_{H^1(\O)}=\int_{\mathop{\rm supp}u_i^n}|\nabla
u_i^n|^2\,dx=\int_{\Omega}\chi_{\{\mathop{\rm supp}u_i^n\}}(x)f_i(x,
u_i^n(x))u_i^n(x)\,dx.
$$
By the choice of $\O^n$, Theorem \ref{limite}(iii), and pointwise
convergence of $u_i^n$ to $u^0_i$, it follows that
$\chi_{\{\mathop{\rm supp}u_i^n\}}\to \chi_{B_i}$ a.e. in $\O$. Hence
the uniform $L^{\infty}$-bound provided in Section
\ref{sec:unif-linfty-bounds} and the Dominated Convergence Theorem
allows to pass to the limit in the right hand side, thus obtaining
$$
\lim_{n\to+\infty}\|u_i^n\|_{H^1(\O)}^2=\int_{B_i}f_i(x,
u_i^0(x))u_i^0(x)\,dx=\|u_i^0\|_{H^1(\O)}^2.
$$
Strong $H^1(\O)$-convergence follows now from weak convergence and
convergence of norms.
\end{rem}

%(ESPANDERE: tutta l'analisi asintotica vale per modelli competitivi)

\section{Coexistence in the Lotka-Volterra models.}\label{perturb}
This section is devoted to prove the existence of solutions to the
auxiliary system when the domain is sufficiently close to $\O^0$ and
the interspecific competition is sufficiently strong.  Precisely, we
shall prove
\begin{teo}\label{t:nostro}
  For any $\vk$ and $n$ sufficiently large, the system with barriers
  \eqref{nostro} admits a solution
  $U^{\vk}=(u_1^{\vk},\dots,u_k^{\vk})\in \big(H^1_0(\Omega^n)\big)^k$
  which is close to $U^0=(u_1^{0},\dots,u_k^{0})$ in
  $\big(H^1_0(\Omega^n)\big)^k$.
\end{teo}

\noindent
In light of Remark \ref{LMpiu}, the above  theorem immediately provides
\begin{coro}
  For any $\vk$ and $n$ sufficiently large, the Lotka-Volterra system
  \eqref{modelLK} admits a solution
  $U^{\vk}=(u_1^{\vk},\dots,u_k^{\vk})\in \big(H^1_0(\Omega^n)\big)^k$
  which is close to $U^0=(u_1^{0},\dots,u_k^{0})$ in
  $\big(H^1_0(\Omega^n)\big)^k$ and satisfies $u_i\geq 0$ for all $i$.
\end{coro}
The proof of Theorem \ref{t:nostro} is obtained by using a standard
topological degree technique (see e.g. \cite{fonsecagangbo}) and it is
based on the ideas introduced in \cite{d} in order to control the
perturbation of the domain.  As a first step, we introduce suitable
operators which allow to reformulate the existence of solutions to
\eqref{problpiu} as a fixed point problem.  For all integers
$n=0,1,\dots,$ we define
\begin{equation*}%\label{eq:ank}
A^{n,\vk}:\ \big(H^1_0(\Omega)\big)^k\to \big(H^1_0(\Omega)\big)^k,\quad
A^{n,\vk}:=L^n\circ F^{n,\vk}\circ i^n,
\end{equation*}
where $i^n: \big(H^1_0(\Omega)\big)^k\! \to\!
\big(H^1(\Omega^n)\big)^k$ is the restriction
$i^n(u_1,\dots,u_k)\!=\!(u_1|_{\O^n},\dots,u_k|_{\O^n})$,
\begin{align*}
  &F^{n,\vk}:\ \big(H^1(\Omega^n)\big)^k \to \big(H^{-1}(\Omega^n)\big)^k, \\
  &F^{n,\vk}(U)=f_i(\cdot,[u_i+u_i^0]^+-u_i^0)\dis- \varkappa
  [u_i+u_i^0]^+\sum_{j\neq i}[u_j+u_j^0]^+,
\end{align*}
and
$$
L^n:\ \big(H^{-1}(\Omega^n)\big)^k\to \big(H^{1}_0(\Omega^n)\big)^k
\hookrightarrow \big(H^{1}_0(\Omega)\big)^k
$$
is defined as: $L^n (h_1,\dots,h_k)=(u_1,\dots, u_k)$ if and only if
$-\Delta u_i=h_i$ in $\Omega^n$, $u_i=0$ on~$\partial \Omega^n$, for
all $i=1,\dots,k$.

With the above notation, it turns out that the solutions of
\eqref{problpiu} in $\O^n$ are in 1-1 correspondence with the fixed
points of $A^{n,\vk}$. We are going to prove the existence
of fixed points of $A^{n,\vk}$ by showing that the
Leray-Schauder degree of the map $\id-A^{n,\vk}$ in a small ball centered
at $U^0$ is different from $0$. We recall that
the Leray-Schauder degree is well defined for operators which differ
from the identity for a compact map. To this aim,
we notice that  it is not restrictive to assume that
$A^{n,\vk}$ is compact from $\big(H^1_0(\Omega)\big)^k$ into itself.
Indeed, if $N<6$, the growth of the nonlinearity $\tilde
q=\max\{2,q\}$ is subcritical, i.e. $\tilde q<\frac{N+2}{N-2}$ and
compactness is guaranteed by the Sobolev-Rellich embedding Theorem.
Otherwise, for $N\geq 6$, using the $L^{\infty}$ bounds proved in Section
\ref{sec:unif-linfty-bounds},
compactness can be recovered by truncating the
coupling term, thus obtaining a subcritical nonlinearity
without affecting the proofs.

The following lemma allows to compute the topological degree of the
unperturbed problem.  In the sequel, we will use the notation $A'(U)$
to denote the Fr\'ech\'et derivative at $U\in X$ of any differentiable
map $A$ from
a Banach space $X$ to the Banach space~$Y$.
\begin{lemma}\label{l:ker}
  Let $\eps>0$ as in assumption \textbf{(ND)}.  There exists $\bar\vk$
  such that for all $\kappa>\bar\vk$, the eigenvalues of
  $\id-\big(A^{0,\vk}\big)'(U^0)$ in $\big(H^1_0(\O^0)\big)^k$ are
  greater than $\eps$.  In particular, the kernel of
  $\id-\big(A^{0,\vk}\big)'(U^0)$ is trivial.
\end{lemma}
\proof
Let us  preliminarly notice that, by Lemma \ref{l:frechet} in
the appendix, the map
$$
F^{0,\vk}:\ \big(H^1(\Omega^0)\big)^k \to \big(H^{-1}(\Omega^0)\big)^k
$$
is Fr\'echet differentiable at $U^0$ and
\begin{equation}\label{eq:10}
\big(F^{0,\vk}\big)'(U^0)[V]=\mathop{\rm Jac}G^{\vk}(U^0)V,\quad
\text{for all }V\in \big(H^1(\Omega^0)\big)^k,
\end{equation}
where
\begin{align*}
G^{\vk}:\ \R^k\to\R^k, \quad
G^{\vk}(U)=f_i(\cdot,u_i)- \varkappa u_i\sum_{j\neq i}u_j
- \varkappa u_i^0\sum_{j\neq i}u_j- \varkappa u_i\sum_{j\neq i}u_j^0,
\end{align*}
and  $ \mathop{\rm Jac}G^{\vk}(U^0)$ denotes the Jacobian matrix of
$G^{\vk}$ at $U^0$.

Let us set ${\mathcal L}_{\vk}:=\id-\big(A^{0,\vk}\big)'(U^0)$ and
write $(H^{1}_0(\Omega^0))^k$ as the direct sum
$$
(H^{1}_0(\Omega^0))^k=\bigoplus_{i=1}^k {\mathcal H}_i,
$$
where
$$
{\mathcal H}_i= H_0^1(\Bi)\times H_0^1(B_{i+1(\textrm{mod }k)})\times
H_0^1(B_{i+2(\textrm{mod }k)})\times\cdots\times
H_0^1(B_{i+k-1(\textrm{mod }k)}).
$$
Spaces ${\mathcal H}_i$ are mutually orthogonal and ${\mathcal
  L}_{\vk}|_{{\mathcal H}_i}:\ {\mathcal H}_i\to {\mathcal H}_i$, so
that it is enough to prove that  $0$ is not an eigenvalue of ${\mathcal
  L}_{\vk}|_{{\mathcal H}_i}$ for all $i=1,\dots,k$.

If $\l$ is an eigenvalue of ${\mathcal L}_{\vk}$ in ${\mathcal H}_1$,
then there exists $V=(v_1,\dots,v_k)\in {\mathcal H}_1$ such that
$(v_1,\dots,v_k)\not=(0,\dots,0)$ and
$$
-(1-\l)\Delta V=\mathop{\rm Jac}G^{\vk}(U^0)V,$$
 i.e.
\begin{equation}\label{eq:5}
-(1-\l)\Delta v_i=\Big(f_i'(\cdot,u_i^0)-2\vk\sum_{j\neq i}u_j^0\Big)v_i-2\vk
u_i^0\sum_{j\neq i}v_j, \quad\text{in }\Omega^0,
\end{equation}
for all $i=1,\dots,k$.
Since $(v_1,\dots,v_k)\not=(0,\dots,0)$, there exists $\ell$ such that
$v_{\ell}\not\equiv 0$. Equation (\ref{eq:5}) for $i=\ell$ in
$B_{\ell}$ reads as
$$
-(1-\l)\Delta v_{\ell}=f_{\ell}'(\cdot,u_{\ell}^0)v_{\ell}, \quad\text{in }B_{\ell},
$$
hence $\l\geq \eps$ in view of assumption {\bf (ND)}.

If $\l$ is an eigenvalue of ${\mathcal L}_{\vk}$ in ${\mathcal H}_i$
for $i\neq1$, then there exists $V=(v_1,\dots,v_k)\in {\mathcal H}_i$,
$V \not=(0,\dots,0)$, which solves (\ref{eq:5}). Let $\ell$ be such
that $v_{\ell}\not\equiv 0$, then equation (\ref{eq:5}) in
$B_{i+\ell-1}$ reads as
$$
-(1-\l)\Delta v_{\ell}=\big(f_{\ell}'(\cdot,0)-2\vk
u^0_{i+\ell-1}\big)v_{\ell},
\quad v_{\ell}\in H^1_0( B_{i+\ell-1}).
$$
Testing the above equation with $v_{\ell}$ we find
\begin{align*}
(1-\l)&\int_{B_{i+\ell-1}}|\nabla
v_{\ell}|^2\,dx=\int_{B_{i+\ell-1}}\big(f_{\ell}'(\cdot,0)-2\vk
u^0_{i+\ell-1}\big)v_{\ell}^2\,dx\\
&\leq
\int_{B_{i+\ell-1}}\big(f_{\ell}'(\cdot,0)-2\vk
u^0_{i+\ell-1}\big)^+v_{\ell}^2\,dx\\
&\leq S^{-1}\bigg( \int_{B_{i+\ell-1}}|\nabla v_{\ell}|^2\,dx \bigg)
\|\big(f_{\ell}'(\cdot,0)-2\vk
u^0_{i+\ell-1}\big)^+\|_{L^{N/2}( B_{i+\ell-1})}
\end{align*}
where $S$ is the best
constant in the Sobolev embedding. Therefore
\begin{equation}\label{eq:6}
\l\geq 1-S^{-1}\|\big(f_{\ell}'(\cdot,0)-2\vk
u^0_{i+\ell-1}\big)^+\|_{L^{N/2}( B_{i+\ell-1})}.
\end{equation}
By the Dominated Convergence Theorem, $\|\big(f_{\ell}'(\cdot,0)-\vk
u^0_{i+\ell-1}\big)^+\|_{L^{N/2}( B_{i+\ell-1})}\to 0$ as
$\vk\to+\infty$ for any $\ell$ and $i$, hence we can find $\bar\vk$ such
that for all $\vk\geq \bar \vk$, for all $i$ and $\ell$
$$
\|\big(f_{\ell}'(\cdot,0)-2\vk
u^0_{i+\ell-1}\big)^+\|_{L^{N/2}( B_{i+\ell-1})}<S(1-\eps).
$$
With this choice of $\bar\vk$, from (\ref{eq:6}) it follows that if
$\l$ is an eigenvalue of ${\mathcal L}_{\vk}$ in ${\mathcal H}_i$ for
$i\neq1$, then $\l\geq\eps$; in particular $\l\neq0$. The proof is
thereby complete.
\endproof
\begin{lemma}\label{l:homot}
There exist $\bar\vk$ and $\bar n$ such that, for all $n\geq \bar n$
and $\vk\geq\bar\vk$,  there holds
$$
U\not=t A^{0,\vk}(U)+(1-t)A^{n,\vk}(U)
$$
for all $t\in[0,1]$ and $U\in (H^1_0(\Omega)\big)^k$ such that
$\|U-U^0\|_{ (H^1_0(\Omega))^k }=\delta$.
\end{lemma}
\proof
Arguing by contradiction, we assume there exist sequences
$n_j\to\infty$ and  $\vk_j\to\infty$, $t_{j}\in[0,1]$, and
$U^j\in (H^1_0(\Omega)\big)^k$ such that $\|U^j-U^0\|_{
  (H^1_0(\Omega))^k}=\delta$ and
\begin{equation}\label{eq:1}
U^j=t_j A^{0,\vk_j}(U^j)+(1-t_j)A^{n_j,\vk_j}(U^j).
\end{equation}
Since $A^{n_j,\vk_j}$ takes values in $\big(H^{1}_0(\Omega^{n_j})\big)^k$, we
have that $U^j\in \big(H^{1}_0(\Omega^{n_j})\big)^k$. Taking the laplacian
of both sides in (\ref{eq:1}), we obtain that $U^j$ solve
\begin{equation*}
\begin{cases}
-\Delta U^j=F^{\vk_j}(U^j), \\[5pt]
U^j\in \big(H^{1}_0(\Omega^{n_j})\big)^k.
\end{cases}
\end{equation*}
Since $\{U^j\}_j$ is bounded in $(H^1_0(\Omega)\big)^k$, up to a
subsequence, $U^j$ converges weakly in $\big(H^{1}_0(\Omega)\big)^k$
to some $U=(u_1,\dots,u_k)\in \big(H^{1}_0(\Omega)\big)^k$.  By
Theorem \ref{limite}, we know that $u_i\cdot u_j=0$ for $i\neq j$,
$u_i\geq 0$ in $\Omega$ and that the $k$-tuple $(u_1,\dots,u_k)$
solves the differential inequality \eqref{sopra}, i.e.
$$
-\Delta \widehat u_i\geq \widehat f_i(x,\widehat u_i), \quad\text{in }\O.
$$
As a matter of fact, arguing as in Theorem \ref{uniqueness} (see the
proof of the Claim 1), it is possible to prove that that $u_i\in
H^{1}_0(\Omega^0)$, hence $U\in\mathcal S(\O^0)$.
Furthermore, since the convergence of $U^j$ to $U$ is actually strong
in $\big(H^{1}_0(\Omega)\big)^k$ by Lemma \ref{limiteforte}, we have
$\sum_i\|u_i-u_i^0\|_{H^1_0(B_i)}^2=
\|U-U^0\|_{(H^1_0(\Omega))^k}^2=\delta^2>0$. This implies the
existence of $i$ such that $u_i\not\equiv u_i^0$, in contradiction
with Theorem \ref{essevuota}.
\endproof

\noindent We now have all the ingredients
to conclude the proof of Theorem \ref{t:nostro}.

\medskip
\noindent
{\bf Proof of Theorem \ref{t:nostro}.}
In view of Theorem \ref{c:isol}, for all $\vk\geq \bar \vk$ we can compute
the Leray-Schauder degree
\begin{equation*}
\deg\big(\id-A^{0,\vk}, B_{(H^1_0(\Omega))^k}(U^0,\delta), 0\big).
\end{equation*}
We learn by Lemma \ref{l:ker} that it turns out to be equal to $+1$.
In light of Lemma~\ref{l:homot}, for $n\geq\bar n$ and $\vk\geq\bar
\vk$, it makes sense to compute the Leray-Schauder degree
$$
I=\deg\big(\id-A^{n,\vk}, B_{(H^1_0(\Omega))^k}(U^0,\delta), 0\big).
$$
By the homotopy invariance property, we have that
$$
\deg\big(\id-A^{n,\vk}, B_{(H^1_0(\Omega))^k}(U^0,\delta), 0\big)=
\deg\big(\id-A^{0,\vk}, B_{(H^1_0(\Omega))^k}(U^0,\delta), 0\big)
$$
and hence $I=+1$.  As a consequence $A^{n,\vk}$ has a fixed point in
$B_{(H^1_0(\Omega))^k}(U^0,\delta)$ which provides a solution
$U=(u_1,...,u_k)$ to~\eqref{problpiu} in~$\Omega^n$, which is close to
$U^0$.  To conclude the proof, it only remains to show that $U$ is a
solution to \eqref{nostro} with the original nonlinearity $f_i$, and
this simply follows from (\ref{pinzo}).
\endproof

\par\noindent
Collecting all the results so far obtained, we can finally
prove or main theorems.

\medskip
\noindent
\textbf{Proof of Theorem 1.}  For a fixed sufficiently large $n$,
let us consider the sequence $U^{\vk}$ of solutions to \eqref{nostro} as in
Theorem \ref{t:nostro}. As $\vk\to\infty$, thanks to Theorem
\ref{limite} we know that $U^{\vk}$ converges strongly to some
$U=(u_1,\dots,u_n)$ in $(H^1_0(\Omega^n))^k$,
such that $U$ is  $H^1$-close to $U^0$, $u_i\geq 0$
for all $i$, $u_i\cdot u_j=0$ if $i\neq j$, $U$ has the non invading
property and satisfies the differential inequalities \eqref{2kvar}.
Hence $U$ belongs to $\mathcal S(\Omega^n)$. The uniqueness is
ensured by Theorem \ref{uniqueness}.
\endproof

\medskip
\noindent
\textbf{Proof of Theorem 2.}  The existence of a solution $U^\kappa$
close to $U^0$ for the two systems is proved in Theorem \ref{t:nostro}
and the subsequent corollary.  The asymptotic analysis as
$\kappa\to\infty$ has been carried out for \eqref{nostro} in Section
\ref{asintotic} and all the results directly come from Theorem
\ref{limite}.  For the Lotka-Volterra model \eqref{modelLK} we still
have that $U^k$ converges to an element of $\mathcal S(\Omega^n)$ by
Remark \ref{LMasynt}.
\endproof

\subsection{Concluding remarks.}
In this paper we have restricted our discussion to homogeneous
Dirichlet boundary conditions.  It has to be stressed that the
technique here employed cannot be used to treat the Neumann no-flux
boundary conditions
\[
\frac{\partial u_i}{\partial \nu}=0,\quad\text{ on }\partial\Omega^n,
\]
the two major obstacles being the difficulty in constructing suitable
extension operators and the lack of continuity of the eigenvalues of
the Laplacian under Neumann boundary conditions with respect to the
perturbation of the domain.  This will be object of forthcoming
studies.

On the other side, our results can be immediately extended to a great
variety of competitive models, not necessarily of Lotka-Volterra type,
since they essentially depend only on the validity of the differential
inequalities \eqref{sopra}.
\\

\noindent\textbf{Acknowledgements.} The authors wish to express their
gratitude to Prof. Susanna Terracini for interesting comments and
discussions. They are also indebted with the anonymous referees for
their helpful remarks.

\section{APPENDIX}
In this appendix, we collect some lemmas used throughout the paper
and the proofs of our  most technical results.
The following simple lemma is needed to prove
Theorem \ref{c:isol}, i.e. to prove isolation of $U^0$.
\begin{lemma}\label{ff}
For all $\eps>0$ there exists $\delta>0$ such that for all
$i=1,\dots,k$ and $u\in
H^1_0(\Omega^0)$,  $\|u-u_i^0\|_{H^1_0(\Omega^0)}\leq \delta$ implies
\begin{equation}
  \label{fder}
  \left|\int_\UB\!\!\big(f_i(x,u_i)-f_i(x,u_i^0)
    -f_i'(x,u^0_i)(u_i-u_i^0)\big)(u_i-u_i^0)\,dx\right|\leq
  \eps\|u_i-u_i^0\|_{H^1_0(\Omega^0)}^2.
\end{equation}
\end{lemma}
\proof
Call $I$ the integral in \eqref{fder}. We can estimate  as follows:
\begin{align*}
  I&=\bigg|\int_{\UB}\bigg[\int_0^1\big(f_i'(x,tu_i+(1-t)u_i^0)
  -f_i'(x,u_i^0)\big)(u_i-u_i^0)^2\,dt\bigg]\,dx\bigg|\\
  &\leq \|u_i-u_i^0\|_{L^{2^*}(\UB)}^2\int_0^1
  \|f_i'(\cdot,tu_i+(1-t)u_i^0)-f_i'(\cdot,u_i^0)\|_{L^{N/2}(\UB)}dt\\
  &\leq S^{-1}\|u_i-u_i^0\|_{H^1_0(\UB)}^2
  \int_0^1\|f_i'(\cdot,tu_i+(1-t)u_i^0)-f_i'(\cdot,u_i^0)\|_{L^{N/2}(\UB)}dt
\end{align*}
where $S$ is the best constant of the Sobolev embedding
$H^1_0\hookrightarrow L^{2^*}$. By continuity of the Nemytskij operator
$f_i':H^1_0(\UB)\to L^{N/2}(\UB)$, $u\mapsto f_i'(\cdot,u(\cdot))$,
there exists $\delta>0$ such that
\begin{align}\label{eq:9}
\|w-u_0^i\|_{H^1_0(\UB)}\leq\delta\Rightarrow\|
f_i'(\cdot,w)-f_i'(\cdot,u_i^0)\|_{L^{N/2}(\UB)}\leq S\eps.
\end{align}
This provides the proof.
\endproof

Let us now prove Theorem \ref{c:isol}, which has played a crucial role
in the degree argument developed in Section \ref{perturb},
 as it ensures the isolation of
the solution to the unperturbed problem.

\medskip
\noindent\textbf{Proof of Theorem \ref{c:isol}.}
Assume that there exists a sequence $U^\varkappa$ of solutions to
\eqref{nostro} such that $u_i^{\vk}> -u_i^0$ for all $i$ and
$U^\varkappa\to U^0$ as $\varkappa\to \infty$.
Let us set $V^\vk=U^\vk-U^0$ and by subtracting the respective
differential equations we obtain, for all $i=1,...,k$:
$$
-\D v^\vk_i=f_i(x,u^\vk_i)-f_i(x,u_i^0)-\,\vk v_i^{\vk}\sum_{j\neq
  i}v_j^{\vk} -2\varkappa u^\vk_i\sum_{j\neq i}u^0_j
-2\varkappa u^0_i
\sum_{j\neq  i}u^\vk_j,
\qquad\text{ in }\Omega^0.
$$
Let us add and subtract the term $f_i'(x,u_i^0)v^\vk_i$; then multiply
by $v_i^\vk$ and integrate on $B_h$ for a fixed $h$.  We have
\begin{eqnarray}
\label{formula}
\int_{B_h}\big[|\nabla
v_i^{\vk}|^2-f_i'(x,u_i^0)|v_i^{\vk}|^2-\big(f_i(x,u_i^{\vk})-f_i(x,u_i^0)-
f_i'(x,u^0_i)v_i^{\vk}\big)v_i^{\vk}\\
\nonumber+
2\varkappa u_i^0(\sum_{j\neq i}v_j^{\vk})v_i^{\vk}+
2\varkappa (\sum_{j\neq i}u_j^0)|v_i^{\vk}|^2+
\varkappa (\sum_{j\neq i}v_j^{\vk})|v_i^{\vk}|^2\big]=0.
\end{eqnarray}
In particular, since  $v_j|_{B_h}=u_j$ if $j\neq h$ while
$v_h|_{B_h}=u_h-u_h^0$, by choosing $h\neq i$ we have
\begin{gather*}
\int_{B_h}|\nabla
v_i^{\vk}|^2-\int_{B_h}\big(f_i(x,u_i^{\vk})-f_i(x,0)-
f_i'(x,0)v_i^{\vk}\big)v_i^{\vk}\\
=\int_{B_h}\big(f_i'(x,0)-\vk u_h^0-\vk\sum_{j\neq i}u_j^{\vk}\big)
|v_i^{\vk}|^2.
\end{gather*}
Let $0<\eps<1$ be given: if $\vk$ is large enough, in light of Lemma
\ref{ff} and since $u_i(x)> 0$ for $x\in B_j$ when $j\neq i$, we know
that
\begin{align*}
  (1-\eps)\int_{B_h}|\nabla v_i^{\vk}|^2&\leq
  \int_{B_h}\big[f_i'(x,0)-\vk
  (u_h^{\vk}+u_h^0)\big]^+ |v_i^{\vk}|^2\\
  &\leq \|[f'_i(x,0)\dis-\varkappa
  (u_h^{\vk}+u_h^0)]^+\|_{L^{N/2}(B_h)}\| v_i^\vk\|^2_{L^{2^*}(B_h)}.
\end{align*}
{\it Claim.
The $L^{N/2}$ norm of $[f'_i(\cdot,0)\dis-\varkappa
(u_h^{\vk}+u_h^0)]^+$  can be made arbitrarily small by letting
$\vk\to\infty$.}

Let us first note that $u_h^{\vk}+u_h^0> 0$ by assumption, hence
$$\|[f'_i(\cdot,0)\dis-\varkappa
(u_h^{\vk}+u_h^0)]^+\|_{L^{N/2}(\omega)}<
\|[f'_i(\cdot,0)]^+\|_{L^{N/2}(\omega)}\leq
\Big(\sup_{x\in\Omega^0}[f_i'(x,0)]^+\Big)
\mu(\omega)^{2/N},
$$
for any measurable $\omega\subset\Omega^0$.
Secondly, since $\|u_h^{\vk}-u_h^0\|_{H^1_0(\Omega^0)}\leq\delta$,
then by the Sobolev embedding
$$
\delta^2\geq\int_{\O^0}|\nabla (u_h^{\vk}-u_h^0)|^2\geq
S\Big(\int_{A_\delta}|u_h^{\vk}-u_h^0|^{2^*}\Big)^{\frac2{2^*}}
\geq \delta \,\mu(A_\delta^{\vk})^{2/2^*}
$$
where
$$
A_\delta^{\vk}=\{x\in B_h:|u_h^{\vk}(x)-u_h^0(x)|^2>\delta\}.
$$
Choose $\delta$ small (independent of $\vk$) enough so
that
$$
\Big(\sup_{x\in\Omega^0}[f_i'(x,0)]^+\Big)
^{N/2}\cdot\mu(A_\delta^{\vk})\leq \frac14\big(S(1-\eps)\big)^{N/2}.
$$

Let us now fix $r>0$ such that
$$
\Big(\sup_{x\in\Omega^0}[f_i'(x,0)]^+\Big)^{N/2}\cdot\mu(B_h\setminus B_h(r))
\leq \frac14\big(S(1-\eps)\big)^{N/2},
$$
where $B_h(r)$ denotes the ball of radius $r$ and with the same center
as $B_h$. We note that there exists
$m>0$ such that $u_h^0(x)\geq m$ for all $x\in  B_h(r)$. Also, for
$0<\sqrt{\delta}<m/2$, we have $u_h^{\vk}+u_h^0>m/2$ in $B_h(r)\setminus
A_\delta^{\vk}$.
With this choice we finally have  $\bar{\vk}$ such that, for all
$\vk\geq\bar{\vk}$, it holds
$[f'_i(x,0)\dis-\varkappa (u_h^{\vk}+u_h^0)]^+(x)=0$ for any $x$ in
$B_h(r)\setminus A_\delta^{\vk}$.
Summing up, the above argument provides
\begin{align*}
  \|[f'_i(\cdot,0)\dis-\varkappa
  (u_h^{\vk}+u_h^0)]^+\|^{N/2}_{L^{N/2}(B_h)} &\leq
  \|[f'_i(\cdot,0)\dis-\varkappa
  (u_h^{\vk}+u_h^0)]^+\|^{N/2}_{L^{N/2}(B_h(r)\setminus A_\delta^{\vk})} \\
  &\quad+\|[f'_i(\cdot,0)\dis-\varkappa
  (u_h^{\vk}+u_h^0)]^+\|^{N/2}_{L^{N/2}((B_h\setminus B_h(r))
    \cup A_\delta^{\vk})}\\
  &\leq \frac12\big(S(1-\eps)\big)^{N/2},
\end{align*}
for $\vk$ large enough, and proves the Claim.
As a consequence, if $\vk$ is large enough, we obtain
$v_i^\vk|_{B_h}\equiv 0$ for all $h\neq i$.  By considering this
information in \eqref{formula} for the choice $h=i$ we get
$$\int_{B_i}\big(|\nabla
v_i^\vk|^2-f_i'(x,u_i^0)|v_i^\vk|^2\big)=
\int_{B_i}\big(f_i(x,u_i^{\vk})-f_i(x,u_i^0)-f_i'(x,u_i^0)v_i^{\vk}\big)v_i^{\vk},
$$
for $\vk$ large enough.  In light of assumption {\bf(ND)} the left
hand side is always bigger than $\eps\|v_i^\vk\|^2$ for some positive
$\eps$. On the other side Lemma \ref{ff} ensures that the right hand
side is lessen by $\eps/2\|v_i^\vk\|^2$ if $\|v_i^\vk\|$ is suitably
small. Hence we reach a contradiction for $\varkappa$ large enough,
unless $v_i^\vk\equiv 0$ for all $i=1,...,k$, that means $U^\vk\equiv
U^0$.
\endproof

\noindent\textbf{Proof of Lemma \ref{sopravuota}.}
Let $i$ be fixed and consider the differential inequality for $\widehat{u}_i$,
$$
-\Delta \widehat u_i\geq \widehat f_i(x,\widehat u_i),\qquad\text{ in } B_i.
$$
Let us test the above inequality with $-\widehat u_i^-$, and denote
$\omega_i:=\{\widehat u_i^->0\}$. This provides
\begin{equation}\label{eq:8}
\int_{\omega_i}|\nabla(\widehat
u_i^-)|^2\leq-\int_{\omega_i}\frac{\widehat f_i(x,\widehat u_i^-)}{\widehat
  u_i^-}(u_i^-)^2\leq M |\mu(\omega_i)|^{2/N}S^{-1}\int_{\omega_i}|\nabla(\widehat
u_i^-)^2|,
\end{equation}
where $M:= \|\widehat f_i(x,\widehat u_i^-)/\widehat
u_i^-\|_{L^\infty}$ is finite by the a priori $L^\infty$ estimate for
$u_i$ as in \eqref{pinzo} and taking into account that $f'_j(0)$
is finite for all $j$. Now, since $\|\widehat
u_i-u_i^0\|_{H^1_0(B_i)}\leq k\delta$,  we have
$$
k^2\delta^2\geq\int_{B_i}|\nabla (\widehat u_i-u_i^0)|^2\geq {\rm
  const\,} \Big(\int_{B_i}|\widehat u_i-u_i^0|^{2^*}\Big)^{2/2^*} \geq
{\rm const\,} \Big(\int_{\omega_i}|u_i^0|^{2^*}\Big)^{2/2^*}.
$$
By absolute continuity of Lebesgue integral, we can choose $\delta$
sufficiently small to ensure that $\mu(\omega_i)<(S/2M)^{N/2}$. Hence
by (\ref{eq:8}) we find
\begin{equation*}
  \int_{\omega_i}|\nabla(\widehat
  u_i^-)|^2\leq \frac12\int_{\omega_i}|\nabla(\widehat
  u_i^-)^2|,
\end{equation*}
which provides $\widehat u_i^-\equiv 0$.
\endproof

\noindent\textbf{Proof of Theorem \ref{essevuota}.}
By Lemma \ref{sopravuota} we know that $\widehat u_i\geq 0$ in
$B_i$ for all $i$. Since $u_j\geq 0$ for all $j$ and the supports are
disjoint, $\sum_{j\neq i}u_j=(\widehat u_i)^-=0$, implying  $\widehat
u_i\equiv u_i$. Hence by coupling the differential inequalities for
$u_i$ and $\widehat u_i$ we obtain that $u_i$ is a solution to
$$-\Delta u_i=f_i(x,u_i),\qquad\text{ in }B_i,$$ with null boundary
conditions. Hence by assumption {\bf (ND)} we obtain $u_i\equiv
u_i^0$.
\endproof

\noindent
The following lemma establishes the Fr\'echet differentiability
of the map $F^{0,\vk}$ defined in Section \ref{perturb}.
\begin{lemma}\label{l:frechet}
  For any $r\in\big[\frac{2N\,\tilde q}{N+2},\frac{2N}{N-2}\big]$, the
  Nemytskij operator
\begin{align*}
  &F^{0,\vk}:\ \big(L^r(\Omega^0)\big)^k \to \big(L^{\frac{r}{\tilde
      q}}(\Omega^0)\big)^k, \\
  &F^{0,\vk}(U)=f_i(\cdot,[u_i+u_i^0]^+-u_i^0)\dis- \varkappa
  [u_i+u_i^0]^+\sum_{j\neq i}[u_j+u_j^0]^+,
\end{align*}
is Fr\'echet differentiable at $U^0$ and
$$
\big(F^{0,\vk}\big)'(U^0)[V]=
\bigg(\Big(f_i'(\cdot,u_i^0)-2\vk\sum_{j\neq i}u_j^0\Big)v_i-2\vk
u_i^0\sum_{j\neq i}v_j\bigg)_{i=1,\dots,k}.
$$
\end{lemma}
\proof
We mean to prove that for all $i=1,\dots,k$
\begin{align*}%\label{eq:3}
  &\Big\| f_i(\cdot,[u_i+u_i^0]^+-u_i^0)\dis-f_i(\cdot,u_i^0)
  -f_i'(\cdot,u_i^0)(u_i-u_i^0)
  \\
  \nonumber &\ -\vk\Big[[u_i+u_i^0]^+\sum_{j\neq
    i}[u_j+u_j^0]^+\!\!-\!  2\Big(\sum_{j\neq i}u_j^0\Big)(u_i-u_i^0)
  \!-\!2 u_i^0\sum_{j\neq i}(u_j-u_j^0)\Big]\Big\|_{(L^{{r}/{\tilde
        q}}(\Omega^0))^k}\\
  \nonumber&=
  o\big(\big\|U-U^0\big\|_{(L^r(\Omega^0))^k}\big)\quad\text{as }
  \big\|U-U^0\big\|_{(L^r(\Omega^0))^k}\to 0.
\end{align*}
We have that
\begin{align*}
  \frac{\Big\| f_i(\cdot,[u_i+u_i^0]^+-u_i^0)\dis-f_i(\cdot,u_i^0)
    -f_i'(\cdot,u_i^0)(u_i-u_i^0)\Big\|^{r/\tilde q}_{(L^{{r}/{\tilde
          q}}(\Omega^0))^k}} {\|u_i-u_i^0\|^{r/\tilde
      q}_{(L^{{r}}(\Omega^0))^k}} \leq I_1+I_2,
\end{align*}
where
\begin{align*}
  &I_1=\frac {\int_{\{u_i+u_i^0>0\}}
    \big|f_i(\cdot,u_i)-f_i(\cdot,u_i^0)
    -f_i'(\cdot,u_i^0)(u_i-u_i^0)|^{r/\tilde q}}{\Big(\int_{\Omega^0}
    \big|u_i-u_i^0|^{r}\Big)^{1/\tilde q}},\\
  &I_2=\frac {\int_{\{u_i+u_i^0<0\}}
    \big|2f_i(\cdot,u_i^0)+f_i'(\cdot,u_i^0)(u_i-u_i^0)|^{r/\tilde
      q}}{\Big(\int_{\{u_i+u_i^0<0\} }
    \big|u_i-u_i^0|^{r}\Big)^{1/\tilde q}}.
\end{align*}
Mimicking the proof of Lemma \ref{ff}, we can easily prove that
$I_1\to0$  as $u_i\to u_i^0$ in $L^r(\Omega)$. Denoting by
$\omega_i:=\{x\in\Omega: u_i(x)+u_i^0(x)<0\}$, we observe that
$$
\int_{\Omega}|u_i-u_i^0|^r\geq\int_{\omega_i}|u_i^0|^r,
$$
hence $|\omega_i|\to 0$ as  $u_i\to u_i^0$ in $L^r(\Omega)$.
From assumptions {\bf (F1)},{\bf (F2)} we have that $f_i(x,s)\leq{\rm
  const\,}(|s|+|s|^{\tilde q})$, hence
\begin{align}\label{eq:4}
  \frac {\Big(\int_{\omega_i} \big|f_i(\cdot,u_i^0)|^{\frac{r}{\tilde
        q}}\Big)^{\frac{\tilde q}r}}{\Big(\int_{\omega_i }
    \big|u_i-u_i^0|^{r}\Big)^{1/ r}}&\leq \frac {\Big(\int_{\omega_i}
    \big|u_i^0|^{\frac{r}{\tilde q}}\Big)^{\frac{\tilde
        q}r}}{\Big(\int_{\omega_i } \big|u_i^0|^{r}\Big)^{\frac1 r}}+
  \frac {\Big(\int_{\omega_i} \big|u_i^0|^{r}\Big)^{\frac{\tilde
        q}r}}{\Big(\int_{\omega_i }
    \big|u_i^0|^{r}\Big)^{\frac1r}}\\
  \nonumber &\leq|\omega_i|^{\frac{\tilde q-1}{r}}+
  \Big(\int_{\omega_i } \big|u_i^0|^{r}\Big)^{\frac{\tilde q-1}{
      r}}=o(1),
\end{align}
as  $u_i\to u_i^0$ in $L^r(\Omega)$. Moreover
\begin{align}\label{eq:7}
  \frac {\Big(\int_{\omega_i}
    \big|f_i'(\cdot,u_i^0)(u_i-u_i^0)|^{\frac{r}{\tilde
        q}}\Big)^{\frac{\tilde q}r}}{\Big(\int_{\omega_i }
    \big|u_i-u_i^0|^{r}\Big)^{1/ r}}&\leq
  \frac{\|u-u_i^0\|_{L^r(\omega_i)}\|f_i'(\cdot,u_i^0)\|_{L^{\frac{r}{\tilde
          q-1}}(\omega_i)}}{\|u-u_i^0\|_{L^r(\omega_i)}}=o(1),
\end{align}
as  $u_i\to u_i^0$ in $L^r(\Omega)$. From (\ref{eq:4}) and
(\ref{eq:7}) it follows that $I_2=o(1)$ as  $u_i\to u_i^0$ in $L^r(\Omega)$.
Hence
\begin{align*}
  \Big\| f_i(\cdot,[u_i+u_i^0]^+\!\!-u_i^0)\!-\!f_i(\cdot,u_i^0) \!
  -\!f_i'(\cdot,u_i^0)(u_i-u_i^0)\Big\|_{(L^{{r}/{\tilde
        q}}(\Omega^0))^k} \!\!=o\big(\|u_i-u_i^0\|_{(L^{{r}/{\tilde
        q}}(\Omega^0))^k}\big)
 \end{align*}
as  $u_i\to u_i^0$ in $L^r(\Omega)$.
On the other hand
\begin{align*}
  &\Big\|[u_i+u_i^0]^+\sum_{j\neq i}[u_j+u_j^0]^+\!\!-\!
  2\Big(\sum_{j\neq i}u_j^0\Big)(u_i-u_i^0) \!-\!2 u_i^0\sum_{j\neq
    i}(u_j-u_j^0)\Big\|_{(L^{{r}/{\tilde q}}(\Omega^0))^k}
  \\
  &\ \leq\sum_{j\neq i}
  \Big(\int_{\Omega^0\setminus(\omega_i\cup\omega_j)}|u_i-u_i^0|^{\frac{r}{\tilde
      q}} |u_j-u_j^0|^{\frac{r}{\tilde
      q}}+\int_{\omega_i\cup\omega_j}|2u_j^0(u_i-u_i^0)+
  2u_i^0(u_j-u_j^0)|^{\frac{r}{\tilde q}}  \Big)^{\frac{\tilde q}r}\\
  &\ \leq {\rm const\,}\sum_{j\neq i} \Big(
  \|u_i-u_i^0\|_{L^r(\Omega)}^{\frac{r}{\tilde q}}
  \|u_j-u_j^0\|_{L^{r/(\tilde q-1)}(\Omega)}^{\frac{r}{\tilde q}}
  \\
  &\qquad\qquad+\|u_i-u_i^0\|_{L^r(\Omega)}^{\frac{r}{\tilde
      q}}\|u_j^0\|_{L^{r/(\tilde
      q-1)}(\omega_i\cup\omega_j)}^{\frac{r}{\tilde q}} +
  \|u_j-u_j^0\|_{L^r(\Omega)}^{\frac{r}{\tilde
      q}}\|u_i^0\|_{L^{r/(\tilde
      q-1)}(\omega_i\cup\omega_j)}^{\frac{r}{\tilde q}}
  \Big)^{\frac{\tilde q}r}
  \\
  &\ =o(1)\quad \text{as } \big\|U-U^0\big\|_{(L^r(\Omega^0))^k}\to 0.
\end{align*}
The proof is thereby complete.
\endproof
Since we are actually working with the truncation $\tilde f_i$ instead
of $f_i$ and $\tilde f_i$ is not $C^1$ with respect to the second
variable, it is worth noticing that this does not create any problem
when linearizing the operator at $U^0$ and the linearization of the
truncated operator is still given by (\ref{eq:10}). Being the proof
very similar to the proof of Lemma~\ref{l:frechet}, we omit it.

\end{document}